\documentclass[10pt]{amsart}

\usepackage{latexsym}
\usepackage{amssymb}
\usepackage{amsmath, amsthm}
\usepackage{amsfonts}

\newtheorem{theorem}{Theorem}

\newtheorem{remark}[theorem]{Remark}

\newtheorem{proposition}[theorem]{Proposition}
\newtheorem{lemma}[theorem]{Lemma}
\newtheorem{corollary}[theorem]{Corollary}

\numberwithin{equation}{section}
\numberwithin{theorem}{section}

\newcommand{\ds}{\displaystyle}

\newcommand{\cR}{\mathbb{R}}

\newcommand{\cC}{\mathbb{C}}

\newcommand{\sfC}{\mathsf{C}}

\newcommand{\cN}{\mathbb{N}}

\newcommand{\cD}{\mathcal{D}}

\newcommand{\ep}{\epsilon}
\newcommand{\de}{\delta}
\newcommand{\al}{\alpha}

\newcommand{\supp}{\mathrm{supp}\,}
\newcommand{\lap}{\triangle}

\begin{document}

\title[Concentration phenomenon for $L^2$-critical NLS]
{A remark on the concentration phenomenon for the $L^2$-critical
nonlinear Schr\"odinger equations}

\author{James Colliander}
\thanks{J.C. is supported in part by N.S.E.R.C. Grant R.G.P.I.N. 250233-03.}
\address{University of Toronto}

\author{Svetlana Roudenko}
\thanks{S.R. is supported in part by N.S.F. grant DMS-0531337.}
\address{Arizona State University}

\begin{abstract}
We observe a link between the window size of mass concentration and
the rate of explosion of the Strichartz norm by revisiting
Bourgain's mass concentration for the $L^2$-critical nonlinear
Schr\"odinger equations.
\end{abstract}

\maketitle

\section{Introduction}

Consider the initial value problem for the $L^2$-critical case of
nonlinear Schr\"odinger equation, $NLS^{\pm}_p(\cR^d)$, $p=\frac4{d}+1$,
\begin{equation}
\label{NLS}
\left\{
\begin{array}{l}
i u_t + \lap u = \pm |u|^{p-1} u,\\
u (0,x) = u_0(x),
\end{array}
\right.
\end{equation}
where $u = u(t,x): \cR \times \cR^d \to \cC$ and
$u_0 \in L^2(\cR^{d})$. $NLS^+_p (\cR^d)$ is called {\it{defocusing}};
$NLS_p^{-} (\cR^d)$ is {\it{focusing}}.

Denote by $[0, T^{\ast})$ the {\it maximal} (forward) existence time
interval of solution $u(t,x)$. For our purposes this means that
any{\footnote{The pair $q = \infty, r=2$ is obviously omitted from
this claim since the $L^2_x$ norm is conserved under the \eqref{NLS}
evolution.}}
Strichartz norm $\ds \Vert u \Vert_{L^q_t
L^r_x{([0,T^{\ast}]\times\cR^d})} = \infty$ and $\ds \| u \|_{L^q_t
L^r_x{([0,t] \times \cR^d)}} < \infty$ for all $t <
T^{\ast}$. Here, the pair $(q,r)$ is 
{\it admissible}, i.e. $\frac2{q}+\frac{d}{r} =\frac{d}{2}$. 

In a breakthrough work \cite{B98}, Bourgain established the mass
concentration phenomenon for finite time blowup solutions of the
cubic $NLS$ in $d=2$ with an $L^2$ initial data (independent of
focusing or defocusing case): consider $NLS^{\pm}_3(\cR^2)$ in
(\ref{NLS}) with $u_0 \in L^2(\cR^2)$; if the blow up time $T^\ast <
\infty$, then $L^2$-norm concentration on a parabolic window occurs
$$
\limsup_{t \nearrow T^{\ast}} \sup_{a \in \cR^2}
\int_{B(a, c(T^{\ast}-t)^{\frac12})} |u(t,x)|^2 \, dx \geq \ep,
$$
where $\ep = \Vert u_0 \Vert^{-M}_{L^2}$ for some $M>0$. The proof used a
refinement of the Strichartz estimate obtained by Moyua, Vargas and Vega
\cite{MVV96} and mass conservation. Compactness properties of blowup
solutions modulo symmetries were obtained by Merle and Vega \cite{MV98}
for the two dimensional case.
Recently, Carles and Keraani in
\cite{CK05} obtained the corresponding results for quintic $NLS$ in $d=1$
and B\'egout and Vargas in \cite{BV05} extended the mass
concentration and compactness modulo symmetries results
for all dimensions $d$ with nonlinearities
$p=\frac4{d}+1$  ($L^2$-critical). The necessary refinement of the
Strichartz inequality for all dimensions comes from the bilinear
Fourier restriction theorem obtained by Tao \cite[Theorem
1.1]{T03}.

In this note we investigate the dependence of the window of mass
concentration upon the growth of the
$L^{\frac{2(d+2)}{d}}_{t,x}([0,t] \times\cR^d)$-norm. We show that
if, close to the blow up time $T^\ast$, the
$L^\frac{2(d+2)}{d}_{t,x}([0,t]\times \cR^d) $-norm grows in time no
slower than ${(T^\ast - t)^{-\beta}}$, then the window of
concentration is of width $(T^\ast - t)^{\frac12+\frac{\beta}2}$. We
also obtain the opposite direction, namely, if the mass
concentration has the concentration window of size $(T^\ast -
t)^{\frac12 + \frac{\beta}2}$, then the growth of
$L^\frac{2(d+2)}{d}_{t,x}([0,t], \times \cR^d) $-norm is no slower
than ${(T^\ast -t)^{-\beta}}$. For the first direction we revisit
the argument of Bourgain and the extension to all dimensions by
B\'egout and Vargas. For the opposite direction we use a restriction
on frequencies (which shows up implicitly in Bourgain's argument) in
order to connect the $L^2_x$-concentration with the space-time
$L^\frac{2(d+2)}{d}_{t,x}$-norm explosion. We also generalize the
above results to the setting of non-polynomial growth and
concentration rates. The result we obtain shows that if the
$L^{\frac{2(d+2)}{d}}_{t,x}$ norm explodes like $f( T^* - t)$ for
certain $f \nearrow \infty$ as $ t \nearrow T^*$, then the
concentration window shrinks at the rate $[-(\partial_t f) (T^* -
t)]^{-\frac{1}{2}}$ and vice versa. As a corollary to \cite{B98} and
\cite{BV05} which proved that parabolic mass concentration occurs,
we obtain that the blow up of the diagonal Strichartz norm must be
at least as fast as $|\ln (T^*-t)|$ (see Corollary \ref{C:log}).

\begin{remark}
Explicit blowup solutions for \eqref{NLS} in the focusing case have
been obtained as the pseudoconformal image of ground and excited
state solitons. These solutions have mass concentration windows
shrinking like $(T^* - t)^{\frac{\beta+1}{2}}$ and their Strichartz
norm $L^{\frac{2(d+2)}{d}}_{t,x} ([0,t] \times \cR^d)$ explodes like
$(T^* - t)^{-\beta}$ with $\beta =1$. Another family of blowup
solutions is known (see \cite{P01} and \cite{MR05}, \cite{MR03})
which concentrates mass
slightly faster (by $~\sqrt{\log |\log (T^* - t)|}$) than
$\beta = 0$. It would be interesting to observe or rule out other blowup
concentration/explosion rates.
\end{remark}

\begin{remark}
It is conjectured that the defocusing problem \eqref{NLS} with the
minus sign is globally well-posed and scatters for all $L^2$ data.
We hope that the results obtained here may be useful in proving that
no concentration occurs in the defocusing problem. For example, in
light of Corollary \ref{C:log}, global well-posedness and scattering
follows if finite time blowup solutions are shown to have
sub-logarithmic Strichartz norm explosion. Also, a result which
rules out very tight concentration windows would imply upper bounds
on the blowup rate of the Strichartz norm. No general upper bounds
on the rate of blowup are known.
\end{remark}

{\sc Notation.} Denote by $l(J)$ the side length of a cube $J
\subset \cR^d$ and $|J|$ its Lebesque measure; $\cD$ is the set of
dyadic cubes in $\cR^d$ with $\ds \tau_k^j = \ds \prod_{i=1}^d
\left[\frac{k_i}{2^j}, \frac{k_i+1}{2^j}\right)$ a dyadic cube, and
when there is no confusion the indices will be dropped $\tau =
\tau_k^j$; for $a \in \cR^d$ and $r>0$ the set $B(a,r) = \{x \in
\cR^d: |x-a|<r\}$ is an open ball of radius $r$.

For a measurable set $E \in \cR^2$, denote by $P_E$ the Fourier
restriction with respect to the $x$-variable: $\widehat{P_E \psi} =
\hat{\psi} \chi_E$. The linear evolution of the Schr\"odinger
equation in (\ref{NLS}) is denoted by $e^{i t \lap}$, i.e.
$$
e^{i t \lap} f(x) = \int_{\cR^d} e^{2\pi i (x\cdot\xi - 2\pi t
|\xi|^2)} \hat{f}(\xi) \, d \xi.
$$

{\sc Acknowledgements:} This project began while the authors
participated in the Fall 2005 Semester on Nonlinear Dispersive Wave
Equations at M.S.R.I. We thank Monica Visan for comments on an
earlier version of this paper.

\section{Strichartz norm explosion $\implies$ tight concentration window}
\label{SectionL2}

First, we show the dependence of the size of mass concentration
window upon the divergence rate of
$L^\frac{2(d+2)}{d}([0,t]\times\cR^d)$-norm.

\begin{proposition}
\label{Prop2} Suppose that $T^{\ast} < \infty$ and
\begin{equation}
\label{L4bound2} \Vert u
\Vert_{L^{\frac{2(d+2)}{d}}({[0,t]\times\cR^d})} \gtrsim
\frac1{(T^{\ast} - t)^{\beta}}~ ~\mbox{for~ some}~\beta >0.
\end{equation}
Then there exists $\ep>0$ such that $\ep = \| u_0
\|_{L^2(\cR^d)}^{-c(d)}$, where\footnote{For example, for $\cR^2$
the argument in (\cite{B98}) gives $c(d) \sim 292$; see the proof
for general $d$.} $c(d)= O(d^4)$,
\begin{equation}
\label{Mass2} \limsup_{t \nearrow T^{\ast}} \sup_{
\begin{array}{c}
\mathrm{cubes}~J \in \cR^d:\\
l(J)< (T^{\ast}-t)^{\frac12 + \frac{\beta}{2}}
\end{array}
} \int_J |u(t,x)|^2 \, dx \geq \ep.
\end{equation}
Furthermore, for any $0 < t < T^{\ast} ~$ there exists a cube
$\tau(t) \subseteq \cR_\xi^d$ of size $\ds l(\tau(t)) \gtrsim
(T^{\ast}-t)^{-(\frac{1}{2} + \frac{\beta}{2})}$ such that
\begin{equation}
\label{Mass2proj} \limsup_{t \nearrow T^{\ast}} \sup_{
\begin{array}{c}
\mathrm{cubes}~J \in \cR^d:\\
l(J)< (T^{\ast}-t)^{\frac12 + \frac{\beta}{2}}
\end{array}
} \int_J |P_{\tau(t)} u(t,x)|^2 \, dx \geq \ep.
\end{equation}
\end{proposition}

Thus, a lower bound on the Strichartz explosion implies tight mass
concentration along a sequence of times. Moreover, the tight
concentration may be frequency localized to the natural scale.

\begin{proof}
We follow \cite{B98} where the mass concentration in the space
dimension $d=2$ is established keeping in mind the generalization to
all space dimensions from \cite{BV05}.

First, recall that a time sequence $\{t_n\} \nearrow T^{\ast}$ is
chosen such that for any $n$
\begin{equation}
\Vert u \Vert_{L^{\frac{2(d+2)}{d}}((t_n, t_{n+1})\times \cR^d)} =
\eta
 \label{StartEqn}
\end{equation}
for some small $\eta$. The decomposition $\ds [0,T^* ) =
\bigcup_{n=0}^\infty [t_n, t_{n+1})$ will play an important role
throughout this paper. If $\eta > 0$ is small enough, then on the
interval $(t_n, t_{n+1})$ the nonlinear part of the evolution
$u(t_n) \mapsto u(t)$ is insignificant compared to the linear flow
$e^{i (t-t_n) \lap} u(t_n)$:
\begin{equation}
 \label{NonlinearEst}
\Vert u - e^{i (t-t_n)\lap} u(t_n)\Vert_{L^{\frac{2(d+2)}{d}}((t_n,
t_{n+1})\times \cR^d)} \leq \Vert u
\Vert^{\frac4{d}+1}_{L^{\frac{2(d+2)}{d}}((t_n, t_{n+1})\times
\cR^d)} = \eta^{\frac4{d}+1},
\end{equation}
and thus,
\begin{equation}
 \label{LinearEst}
\Vert e^{i (t-t_n)\lap} u(t_n)\Vert_{L^{\frac{2(d+2)}{d}}((t_n,
t_{n+1})\times \cR^d)} \thicksim \eta.
\end{equation}

We impose $\eta < \min(1, \frac1{2c} (T^{\ast})^{-\beta})$, where
$c$ is the implicit constant in (\ref{L4bound2}). Then the bound
(\ref{L4bound2}) implies
\begin{equation}
\label{LowBound} \eta = \Vert u \Vert_{L^{\frac{2(d+2)}{d}}((t_n,
t_{n+1})\times \cR^d)} \gtrsim \frac{(t_{n+1}-
t_n)}{(T^{\ast}-t_n)^{\beta+1}},
\end{equation}
i.e., the sequence $\{t_n\}$ has the following property
\begin{equation}
\label{t-dependence} t_{n+1}-t_n \lesssim \eta \, (T^{\ast} -
t_n)^{\beta+1}.
\end{equation}

Fix $n \in \cN$ and the time interval $(t_n,
t_{n+1})$. Denote $f(x)=u(t_n,x)$, note that $\Vert f
\Vert_{L^2(\cR^d)} = \Vert u_0 \Vert_{L^2(\cR^d)}$ by mass
conservation. Using the Squares Lemma (\S 2 in \cite{B98} and
\cite[Lemma 3.1]{BV05}), we obtain the following localizations
in frequency.
\begin{itemize}
\item
For any $\ep_0 >0$
there exist $N_0 = N_0(\Vert f \Vert_{L^2}, d, \ep_0)$ and a finite
collection $\{f_{j} \}_{j=1}^{N_0} \in L^2(\cR^d)$ with $\supp
\hat{f}_{j} \subseteq \tau_j$ - a cube in $\cR^d$, $l(\tau_j)
\leq c(\Vert f \Vert_{L^2}, \eta, \ep_0)\cdot A_{j}$, $|\hat{f}_{j}|
< {A_j}^{-d/2}$  such that
\begin{equation}
 \label{loc1}
\Vert e^{i t \lap} f - \sum_{j=1}^{N_0} e^{i t \lap}
f_{j}\Vert_{L^{\frac{2(d+2)}{d}}(\cR\times\cR^d)} < \ep_0.
\end{equation}
\end{itemize}

Expand (\ref{StartEqn}) and apply (\ref{NonlinearEst}) and
(\ref{LinearEst})
$$
\eta^{\frac{2(d+2)}{d}} = \int_{(t_n, t_{n+1})\times \cR^d} \left[
u(t,x) \left(\overline{e^{i(t-t_n)\lap} f + (u(t,x) -
e^{i(t-t_n)\lap} f)}\right) \right.
$$
$$
\times \left. \left|e^{i(t-t_n)\lap} f + (u(t,x) - e^{i(t-t_n)\lap}
f)\right|^{\frac4{d}} \right] \, dx \, dt
$$
$$
= \int_{(t_n, t_{n+1})\times \cR^d} u(t,x) \,
(\overline{e^{i(t-t_n)\lap} f}) \, |e^{i(t-t_n)\lap} f|^{\frac4{d}}
\, \,  dx \, dt
$$
$$
+ \, O (\eta^{2(\frac4{d}+1)}) \qquad \text{if} \quad d \leq 4,
$$
or
$$
+ \, O (\eta^{2 +\frac4{d}(\frac4{d}+1)}) \quad \text{if} \quad d
>4.
$$

Choose $\ep_0 = \eta^{\frac{d+4}{d}}$. Then the above estimate
together with (\ref{loc1}) may be rewritten as
\begin{equation}
 \label{Reduction1}
\eta^{\frac{2(d+2)}{d}} =  \!\!\!\int\limits_{(t_n, t_{n+1})\times
\cR^d} \hspace{-.6cm}u(t,x) \, \left(\overline{\sum_{j=1}^{N_0} e^{i
(t-t_n) \lap} f_{j}(x)} \right) \left|\sum_{j=1}^{N_0}
e^{i(t-t_n)\lap} f_j (x) \right|^{\frac4{d}} dx \, dt
\end{equation}
$$
+ \, O (\eta^{2(\frac4{d}+1)}) \qquad \text{if} \quad d \leq 4,
$$
or
$$
+ \, O (\eta^{2 +\frac4{d}(\frac4{d}+1)}) \quad \text{if} \quad d
>4.
$$

Note that ~$\, 2(\frac4{d}+1) > \frac{2(d+2)}{d}\,$~ or ~$\, 2
+\frac4{d}(\frac4{d}+1) > \frac{2(d+2)}{d}\,$~ for any $d>0$. Since
$N_0$ is finite, it follows from (\ref{Reduction1}) that there
exists $j_0$ ($1 \leq j_0 \leq N_0$) such that
\begin{equation}
 \label{Est2}
\frac{1}{2} \frac{\eta^{\frac{2(d+2)}{d}}}{N_0^{1+4/d}} \, \leq
\left|\, \int\limits_{(t_n, t_{n+1})\times \cR^d} \hspace{-.6cm}
u(t,x) (\overline{e^{i (t-t_n) \lap} f_{j_0}(x)}) \,
|e^{i(t-t_n)\lap} f_{j_0}(x)|^{\frac4{d}} \, dx \, dt \right|.
\end{equation}
Denote $\tau = \tau_{j_0}$, $A=A_{j_0}$ $(l(\tau)\leq c_0 \,A)$ and
the center of $\tau$ by $\xi_0$. Since $|\hat{f}_{j_0}| <
{A}^{-d/2}$ and $\hat{f}_{j_0}$ is supported in $\tau$ with $l(\tau)
= c_0(\Vert f \Vert_{L^2}, \eta)\cdot A$, we obtain
$$
\ds \sup_{t} \Vert e^{i (t - t_n) \lap} f_{j_0} \Vert_{L^{\infty}_x}
\leq A^{-d/2}\, |\tau| = c \, A^{d/2} \quad \text{with} \quad
c=c(\Vert f \Vert_{L^2(\cR^d)}, \eta).
$$
Hence, (\ref{Est2}) becomes
$$
\frac12 \frac{\eta^{\frac{2(d+2)}{d}}}{N_0^{1+4/d}} \, \leq A^2 \,
|\hspace{-.6cm}\int\limits_{(t_n, t_{n+1}) \times \cR^d} \!\!\!
u(t,x) (\overline{e^{i (t-t_n) \lap} f_{j_0}(x)}) \, dx \, dt \,| .
$$
Using Plancherel and $\supp \hat{f}_{j_0} \subseteq \tau$, we obtain
a refined version of the previous inequality
\begin{equation}
 \label{RefinedProjection}
\frac{\eta^{\frac{2(d+2)}{d}}}{2 \,N_0^{1+4/d}} \, \leq A^2
\int\limits_{(t_n, t_{n+1}) \times \cR^d} |P_{\tau} u(t,x) | \,
|e^{i (t-t_n) \lap} f_{j_0}(x)| \, dx \, dt.
\end{equation}

Using the  Tubes Lemma (\S 3 in \cite{B98} and \cite[Lemma
3.3]{BV05}), we obtain a further space-time localization:
\begin{itemize}
\item
Let $\ep_1 = (\frac1{4 N_0})^{d/4}\,\eta$. Then there exist
$N_1=N_1(\Vert f \Vert_{L^2}, d, \ep_1)$ and a sequence of tubes
$\{Q_{k}\} = \{I_k \times K_k(t)\} \subseteq \cR \times \cR^d$,
where $|I_k|=\frac1{A^2}$ and $K_k(t) = 4\pi t \xi_0 + \sfC$, $\sfC
\in \cD$, with $l(\sfC) = \frac1{A}$ such that
\begin{equation}
 \label{loc2}
\Vert e^{i (t-t_n) \lap}
f_{j_0}\Vert_{L^{\frac{2(d+2)}{d}}(\cR\times\cR^d \setminus
{\underset{k=1}{\overset{N_1}{\cup}} Q_{k}})} < \ep_1.
\end{equation}
\end{itemize}
By (\ref{RefinedProjection}) and (\ref{loc2}), we obtain
\begin{equation}
 \label{Reduction2}
\frac14 \frac{\eta^{\frac{2(d+2)}{d}}}{N_0^{1+4/d}} \, \leq A^2
\hspace{-.7cm}\int\limits_{(t_n, t_{n+1}) \times \cR^d \cap
({\underset{k=1}{\overset{N_1}{\cup}} Q_k})} \hspace{-.7cm} |P_\tau
u(t,x)|\, |{e^{i (t-t_n) \lap} f_{j_0}(x)}| \, dx \, dt.
\end{equation}
Since the number of tubes $N_1$ is finite, there exists a tube
$Q_{k_0} = I \times K(t)$ such that (\ref{Reduction2}) produces
\begin{equation}
 \label{Est1}
\frac14\frac{\eta^{\frac{2(d+2)}{d}}}{N_0^{1+4/d} \, N_1} \, \leq
A^2 \hspace{-.7cm}\int\limits_{(t_n, t_{n+1}) \times \cR^d \cap
Q_{k_0}} \hspace{-.7cm} |P_\tau u(t,x)| \, |e^{i (t-t_n) \lap}
f_{j_0}(x)| \, dx \, dt .
\end{equation}
Applying Cauchy-Schwarz we obtain
\begin{equation}
 \label{Est3}
\frac{\eta^{\frac{2(d+2)}{d}}}{4 \,N_0^{1+4/d} \, N_1} \, \leq A
\left(\, \int\limits_{(t_n, t_{n+1}) \times \cR^d \cap Q_{k_0}} \!\!
|P_{\tau} u(t,x) |^2 \, dx \, dt \right)^{1/2}.
\end{equation}
Since\footnote{We estimate $c(d) \approx
\frac{(d+4)(d+3)(d^2+3d+4)^2}{(d+1)(d+2)}$.} $N_0^{\frac{d+4}{d}}
N_1 \sim \Vert f \Vert^{c(d)}_{L^2}$, using the conservation of
mass, we get
$$
c \leq A \, \Vert u_0 \Vert_{L^2} \, (t_{n+1} - t_n)^{1/2}, \quad
\text{where} \quad c=c(\Vert u_0 \Vert_{L^2}, d, \eta).
$$
By (\ref{t-dependence}), we obtain
$$
\frac1{A} \leq c \, (T^{\ast} - t_n)^{\frac{\beta+1}2} \quad
\mbox{again with} \quad c = c(\Vert u_0 \Vert_{L^2}, d, \eta).
$$
Considering (\ref{Est3}) again, we have
\begin{align}
 \label{E:conc-sup1}
c \, \eta^{\frac{4(d+2)}{d}} &\leq  A^2 \int_{I \cap (t_n, t_{n+1})}
\int_{K(t)} |P_{\tau} u(t,x)|^2 dx \, dt\\
& \leq \sup_{t \in I\cap(t_n, t_{n+1})} \int_{K(t)} |P_{\tau}
u(t,x)|^2 \, dx.
 \label{E:conc-sup2}
\end{align}
Therefore, there exists a mass concentration time
$t_n^* \in I\cap(t_n, t_{n+1})$ such that
\begin{equation}
 \label{E:conc-sup3}
\int_{K(t_n^*)} |P_{\tau} u(t_n^*, x)|^2 \, dx \geq c'
\eta^{\frac{4(d+2)}{d}}.
\end{equation}
The $limsup$ claim in \eqref{Mass2proj} will be realized along the
sequence $t_n^* \nearrow T^*$. Recall that $l(K(t_n^*)) = \frac1{A}
< c \, (T^{\ast} - t_n)^{\frac{\beta+1}{2}}$, and therefore,
$K(t_n^*) \subseteq B(a, \sqrt d\, c \, (T^{\ast} -
t_n)^{\frac{\beta+1}{2}})$ for some $a \in \cR^d$. Observe that
\begin{equation}
\label{time-0} T^{\ast} - t_n^* > T^{\ast} - t_{n+1} = T^{\ast} -
t_n - (t_{n+1}-t_n)
\end{equation}
$$
\geq (T^{\ast} - t_n)(1 - c \, \eta (T^{\ast} - t_n)^{\beta})
> \frac12 (T^{\ast} - t_n),
$$
where the last estimate follows from $c \,\eta \,(T^{\ast} -
t_n)^{\beta} < c \,\eta \,(T^{\ast})^{-\beta}<\frac12$. Hence, $B(a,
\sqrt d\,c\, (T^{\ast} - t_n)^{\frac{\beta+1}{2}})$ (and thus,
$K(t_n^*)$) can be covered by a finite number of balls (or cubes) of
radius (side length) $(T^{\ast} - t_n^*)^{\frac{\beta+1}{2}}$ (and
this number is independent of $n$). Choosing one of them, and noting
that $\ds l(\tau) \geq \frac{c}{(T^{\ast} -
t_n)^{\frac{\beta+1}{2}}} \approx \frac{c}{(T^{\ast} -
t_n^*)^{\frac{\beta+1}{2}}}$, we get
$$
\int_{B\left(a, \,c(T^{\ast} - t_n^*)^{\frac{\beta+1}{2}}\right)}
|P_{\tau(t_n^*)} u(t_n^*, x)|^2 \, dx
\geq \ep,
$$
and since $n$ is arbitrary, the proposition follows.
\end{proof}

\begin{remark} Observe that we did not use the splitting of the
interval $(t_n, t_{n+1})$ as on page 261 in \cite{B98}, since we had
the estimate (\ref{t-dependence}) of $(t_{n+1}-t_n)$ in terms of
$(T^{\ast}-t_n)^{\beta +1}$, $\beta +1 > 1$, which gives a nonzero
bound in (\ref{time-0}). In Bourgain's argument $\beta = 0$, i.e.
$(t_{n+1}-t_n)< (T^{\ast}-t_n)$, which is not enough to conclude
mass concentration with the above argument, and thus, a more careful
splitting of the time interval is needed.
\end{remark}

Note that the construction of $t_n^*$ given in the proof above
provides more information about the mass concentration than is
claimed in \eqref{Mass2proj}. For example, we know that there is a
concentration time $t_n^*$ in each of the time intervals $[t_n ,
t_{n+1})$. The next statement contains a strengthened conclusion
which shows that the concentration actually holds on a thickened
interval of times containing $t_n^*$ of size proportional to
$t_{n+1} - t_n$.

\begin{corollary}
 \label{CorThickTime}
Assume the hypotheses of Proposition \ref{Prop2}. The conclusion
\eqref{Mass2proj} may be strengthened as follows: There exist ~$0 <
\sigma < 1$\footnote{If $\sigma = 0$, then this statement coincides
with the theorem by choosing $I_n = \{t_n^*\}$.} and a sequence of
time intervals $\{I_n\}$ with $I_n \subset (t_n,t_{n+1})$  and
$|I_n| = \sigma \, (t_{n+1} - t_n)$, uniform for all $n$, such that
for some $\tilde\sigma = \tilde\sigma(\sigma)
>0$\footnote{with $\tilde\sigma = 0$ if $\sigma = 0$ and $\tilde\sigma =
1$ when $\sigma = 1$.}
we have
\begin{equation}
 \label{refinedMass2proj}
\lim_{n \to \infty} \inf_{~t \in I_n} \sup_{ \begin{array}{c}
\mathrm{cubes}~J \in \cR^d:\\
l(J)< (T^{\ast}-t)^{\frac12 + \frac{\beta}{2}}
\end{array}
} \int\limits_{J} |P_{\tau(t)} u(x, t)|^2 \, dx \geq
(1-\tilde\sigma)\, {\ep}.
\end{equation}
\end{corollary}

\begin{proof}
Recall the inequalities (\ref{E:conc-sup1}) - (\ref{E:conc-sup2})
from the proof of Proposition \ref{Prop2}. From (\ref{E:conc-sup2})
only one concentration time $t_n^* \in I \cap (t_n, t_{n+1})$ was
selected such that (\ref{E:conc-sup3}) holds. However,
(\ref{E:conc-sup1}) contains a stronger statement, namely, on each
set $I \cap (t_n, t_{n+1})$ there exists a subset $E_n$ such that
for any $t \in E_n$ we have
\begin{equation}
 \label{E:setEn}
c' \, \eta^{\frac{4(d+2)}{d}} \leq \int_{K(t)} |P_{\tau} u(x, t)|^2
\, dx.
\end{equation}
Note that since $u_0 \in L^2_x$, by the local well-posedness and
mass conservation, $u \in C^0_t(L^2_x)$, and so $P_{\tau} u(t)$ is
also continuous in time, which means that the set $E_n$ above can be
chosen to be an interval, denote it by $I_n$. Next we estimate how
large $I_n$ can be in comparison with $I\cap(t_n,t_{n+1})$. First,
recall that $|I| = \frac1{A^2}$ and $\frac1{A^2} \lesssim (t_{n+1} -
t_n)$. Since the cube $\tau \subset \cR^d$ has the center $\xi_0$
and side length $l(\tau)\leq c_0 A$, the function $P_{\tau} u$ (on
the time interval $(t_n,t_{n+1})$) contains frequencies $\xi \in
\tau$, and thus, $|\xi| \leq c(\xi_0,c_0)\cdot A$. By the
uncertainty principle (for example, p. 332 of \cite{T05}) $P_{\tau}
u$ is approximately constant on spatial balls of radius
$\frac{c}{A}$ for some small $c$, in particular, since $l(K(t)) =
\frac1{A}$, it will be approximately constant on some fixed part of
$K(t)$. By the propagation of Schrodinger waves, this set will
persist for an interval of times of measure $\sim \frac1{A^2}$
(after which it may disperse). This length scale is exactly
comparable with the size of $I$ (note independently of the step
$n$), so we can find $0<\sigma<1$ such that $|I_n| =\sigma\,
|I\cap(t_n,t_{n+1})|$ for all $n$ and (\ref{E:setEn}) holds for all
$t \in I_n$ and some $\tilde\sigma >0$:
\begin{equation}
 \label{E:setIn}
(1 - \tilde\sigma)\, c \, \eta^{\frac{4(d+2)}{d}} \leq \int_{K(t)}
|P_{\tau} u(x, t)|^2 \, dx.
\end{equation}
For the above heuristics we need the following lemma

\begin{lemma}
Let $f \in L_x^2(\cR^d)$ and $\supp \hat{f} \subset [0,1]^d$.
Suppose that for some constant $c_1 > 0$
$$
\int_{[0,1]^d} |f(x)|^2 \, dx \geq c_1.
$$
Then for $|t| < c(c_1, \| f \|_{L^2})$ the same concentration holds
for the linear Schr\"odinger evolution of $f$, i.e.,
$$
\int_{[0,1]^d} |e^{it\lap} f(x)|^2 \, dx \geq \frac{c_1}2.
$$
\end{lemma}

\begin{proof}
A basic calculation yields
\begin{equation}
 \label{E:initialconc}
c_1 \leq \int_{[0,1]^d} |f(x)|^2 \, dx \leq \int_{[0,1]^d}
|\left(f(x) - e^{it \lap} f(x)\right) + e^{it \lap} f(x)|^2 \, dx
\end{equation}
$$
\leq 2 \, \left(\int_{[0,1]^d} |f(x) - e^{it \lap} f(x)|^2 \, dx +
\int_{[0,1]^d} |e^{it \lap} f(x)|^2 \, dx \right) := A + B.
$$
We reexpress the integrand in A using the Fourier transform
$$
|f(x) - e^{it \lap} f(x)|^2 \leq \int_{[0,1]^d} |(e^{-4\pi^2 i t
|\xi|^2} - 1) e^{2 \pi i x \xi} \hat{f}(\xi)|^2 \, d\xi
$$
$$
\leq \sup_{\xi \in [0,1]^d} \left|e^{-4\pi i t |\xi|^2} - 1
\right|^2 \int_{[0,1]^d} |\hat{f}(\xi)|^2 \, d\xi \leq 2 \,(4\pi^2
\, t)^2 \, \Vert f \Vert_{L^2}^2. $$ Here we used the estimate
$$ |e^{4\pi^2 i t} - 1|^2 \leq (\cos (4\pi^2 t) -1)^2 + \sin^2(4\pi^2 t)
\leq 2\, (1-\cos^2(4\pi^2 \, t)) \leq 2\,(4\pi^2 \, t)^2.
$$
If we restrict  $t$ such that $\ds |t| \leq \frac{1}{4 \pi^2 \| f
\|_{L^2}} \sqrt{ \frac{c_1}{8}}$, then $\ds A \leq \frac{c_1}{2}$,
and we obtain the conclusion of the lemma.
\end{proof}

We return to the proof of Corollary \ref{CorThickTime}. The
preceding lemma shows that $L^2$ functions which are band limited to
unit scale and lower frequencies and which are mass concentrated at
unit scale remain mass concentrated at unit scale for unit time
under the linear Schr\"odinger flow. Applying the dilation
invariance shows that $L^2$ functions which are band limited to
frequencies $|\xi| \lesssim A$ and which are mass concentrated on
$|x| \lesssim \frac{1}{A}$ will remain mass concentrated for time
$|t| \lesssim \frac{1}{A^2}$ under the linear Schr\"odinger flow.
Finally, using the translation and Galilean invariances, we observe
that this parabolic mass concentration persistence property holds
without special reference to the frequency or spatial origin.

Now the rest of the argument in the proof repeats for any $t \in
I_n$ (for example, (\ref{time-0}) holds for any $t \in I_n$ because
of the fixed proportion $\sigma$ to $(T^*-t_n)^{\beta+1}$ and we
obtain (\ref{refinedMass2proj}). This completes the proof of
Corollary \ref{CorThickTime}.
\end{proof}

\begin{corollary} \label{GeneralCase}
The above results can be extended to a more general form of the
lower bound on the Strichartz norm in (\ref{L4bound2}). Suppose
$$
\Vert u \Vert_{L^{\frac{2(d+2)}{d}}({[0,t]\times\cR^d})} \gtrsim
G(T^*-t),
$$
where $G(s) \to +\infty$ as $s \to 0$ and $G \in C^1(0,1)$. Then the
window in the mass concentration (\ref{Mass2}) changes as follows:
\begin{enumerate}
\item
if $G(T^*-t) \gtrsim |\ln (T^*-t)|^{\gamma}$ with $\gamma \geq 1$,
then $l(J) < [-(\partial_t G)(T^* - t)]^{-1/2}$,

\item
otherwise, $l(J) < (T^* - t)^{1/2}$.
\end{enumerate}
Similar changes hold for (\ref{Mass2proj}) and also in Corollary
\ref{CorThickTime}.
\end{corollary}

Observe that by the argument of Bourgain \cite{B98} and
B\'egout-Vargas \cite{BV05} we always have the case {\rm (2)}. The
case {\rm (1)} is an improvement of {\rm (2)} when $G$ grows faster
than $|\ln (T^*-t)|$, otherwise, the argument of Proposition
\ref{Prop2} gives a weaker statement, i.e. the width of the window
of concentration given by $G_t$ is narrower than the parabolic window.

\begin{proof}
The proof of this corollary follows the proof of Proposition
\ref{Prop2} (and Corollary \ref{CorThickTime}) with the following
changes. Given $G$ as above, the estimate (\ref{LowBound}) changes
to
$$
\eta = \Vert u \Vert_{L^{\frac{2(d+2)}{d}}((t_n, t_{n+1})\times
\cR^d)} \gtrsim (t_{n+1}- t_n) \, [-(\partial_t G)(T^{\ast}-t_n)],
$$
with the constant independent of $n$, and thus,
$$
t_{n+1}-t_n \lesssim \eta \, [-(\partial_t G)(T^{\ast} - t_n)]^{-1}.
$$
Hence, the size of $\tau = \tau_{j_0}$ is estimated as
$$
\frac1{A} \leq c(\|u_0 \|^2_{L^2}, \eta, d) \, [-(\partial_t
G)(T^*-t_n)]^{-1/2},
$$
which implies the result in {\rm (1)}. Note that if $G$ has a faster
grows than $|\ln(T^*-t)|$, then to get the parabolic window of
concentration, we need the extra splitting of the interval $(t_n,
t_{n+1})$ as on page 261 in \cite{B98} or in Step 3 of Prop. 4.1 in
\cite{BV05}. This finishes the proof.
\end{proof}

As an example, consider $G(T^*-t) = |\ln (T^*-t)|^{1+\ep}$, $\ep>0$,
then $l(J) < (T^* - t)^{1/2} |\ln(T^*-t)|^{-\ep/2}$ which is wider
than the parabolic window. If $G(T^*-t) = \ln|\ln(T^*-t)|$, then the
case (2) holds and the window of concentration is parabolic.

\section{Tight concentration window $\implies$ Strichartz norm explosion}

The following statement shows how the radius of mass concentration
affects the divergence rate of the
$L^{\frac{2(d+2)}{d}}_{t,x}$-norm. We will use the shorthand
notation $P_{L(t)}$ to denote the Fourier restriction operator
$P_{\{|\xi| \leq L(t)\}}$ and $F(t) = \Vert u
\Vert^{\frac{2(d+2)}{d}}_{L^{\frac{2(d+2)}{d}}{([0,t])\times
\cR^d})}$.

\begin{proposition}[Local Estimate]
 \label{propSup}
Let  $u \in C([0,T^{\ast}) ; L^2(\cR^d)) \, \cap \,
L^\frac{2(d+2)}{d}([0,T^{\ast})$ ; $L^\frac{2(d+2)}{d}(\cR^d))$~ be
the maximal solution of $NLS^\pm_p(\cR^d)$, $p=\frac4{d}+1$, with
$u_0 \in L^2(\cR^d)$. \\
For $\ep>0$ let $\ds \kappa(\ep) = 2^{-(d+2)} \left({\ep} \,\Vert
u_0 \Vert^{-2}_{L^2}/ 8 \right)^{1/d}$ and for $\alpha>0$ define
$\ds L(t) = \frac12 \frac{\kappa(\ep)}{(T^{\ast}-t)^{\alpha}}$.

Suppose there exists $\alpha \geq \frac12$ and $\ep
> 0$ such that
\begin{equation}
\label{Mass1Sup}
\limsup_{t \nearrow T^{\ast}} \sup_{
\begin{array}{c}
\mathrm{cubes}~J \subset \cR^d:\\
l(J)< (T^{\ast}-t)^\alpha
\end{array}
} \int_J |P_{L(t)} u(t,x)|^2 \, dx \geq \ep.
\end{equation}

Then there exists $t_n \nearrow T^{\ast}$ such that
\begin{equation}
\label{L4boundSup} F'(t_n) \gtrsim \frac{1}{(T^*-t_n)^{2\alpha}}.
\end{equation}

\end{proposition}

\begin{remark}
If in (\ref{Mass1Sup}) one has $\liminf$ instead of $\limsup$, then
(\ref{L4boundSup}) holds for any sequence $t_n \nearrow T^\ast$.
\end{remark}

\begin{proof}
The condition (\ref{Mass1Sup}) implies that there exists a sequence
of times $\{t_n\}_{n=1}^{\infty}$ with $t_n \nearrow T^{\ast}$ and a
sequence of cubes $\{J_n\}_{n=1}^{\infty} \subset \cR^d$ with
$l(J_n) < (T^{\ast} - t_n)^\alpha$ such that
\begin{equation}
 \label{limsup1}
\frac{\ep}{2} \leq \int_{J_n} |P_{L(t_n)} u(x,t_n)|^2 \, dx.
\end{equation}
Since $F(t)$ is an increasing function, by the monotonicity theorem
(e.g. see \cite{R}), it follows that $F'$ exists for a.e. t. We may
assume that $F'(t_n)$ exists and is finite: for any $\tilde\ep
> 0$ there exists $t_n^* \in (t_n, t_n+\tilde\ep)$ such that
$F'(t_n^*)$ exists and is finite, so we choose $\tilde\ep_{n}$ such
that $\ds \tilde\ep_n < 2^{d(d+3)-2} \, \ep \, \frac{(T^* -
t_n)^{d\, \alpha}}{(T^* - (t_n + \ep))^{d \alpha}}$ (note that
$\tilde\ep_n \to 0$ as $n \to \infty$). Then using the triangle
inequality in (\ref{limsup1}), we obtain
$$
\frac{\ep}2 \leq 2\, \left( \int\limits_{~J_n} |P_{L(t_n)}u(x,t_n) -
P_{L(t_n^*)} u(x,t_n^*)|^{2} \, dx + \int\limits_{J_n} |P_{L(t_n^*)}
u(x,t_n^*)|^{2}\, dx \right)
$$
$$
\leq \frac{\ep}4 + 2 \, \int\limits_{J_n} |P_{L(t_n^*)}
u(x,t_n^*)|^{2}\, dx,
$$
where the bound on the first term in the right hand side is
discussed below, and thus, we obtain (\ref{limsup1}) with
$\frac{\ep}4$ on the left-hand side. Observe that
\begin{equation}
 \label{E:difference}
|P_{L(t_n)}u(x,t_n) - P_{L(t_n^*)} u(x,t^*_n)| \leq \int_{|\xi|\leq
L(t_n^*)} |\hat{u}(\xi, t_n) - \hat{u}(\xi, t_n^*)| \, d\xi
\end{equation}
$$
\leq \left( \int_{\cR^d} |\hat{u}(\xi, t_n) - \hat{u}(\xi, t_n^*)|^2
\, d\xi \right)^{1/2} \left(\int_{\{\xi \in \cR^d: |\xi|\leq \,
L(t_n^*)\}} 1 \, d\xi \right)^{1/2}
$$
$$
\leq 2 \, \Vert u_0 \Vert_{L^2(\cR^d)} \, [2 L(t_n^*)]^{d/2} \leq
\left( \frac{\tilde\ep_n}2 \, \frac{1}{2^{d(d+3)}} \,
\frac{1}{(T^*-t_n^*)^{d \, \alpha}} \right)^{1/2},
$$

and integrating over $J_n$, we obtain
$$
2 \int\limits_{~J_n} |P_{L(t_n)}u(x,t_n) - P_{L(t_n^*)}
u(x,t_n^*)|^{2} \, dx
 \leq \frac{\tilde\ep_n}{2^{d(d+3)}} \,\frac{{(T^* -
 t_n)^{d
\alpha}}}{(T^* - t_n^*)^{d\, \alpha}} < \frac{\ep}4,
$$
by the choice of $\tilde\ep_n$. Now we may re-denote the sequence
$t_n^*$ by $t_n$.

Returning to (\ref{limsup1}), we obtain
\begin{equation}
 \label{limsup2}
\frac{\ep}2 \leq \left( \int_{J_n}|P_{L(t_n)}
u(x,t_n)|^{\frac{2(d+2)}{d}} \, dx \right)^{\frac{d}{d+2}}
l(J_n)^{\frac{2d}{d+2}},
\end{equation}
by H\"older's inequality.

Fix $n \in \cN$ and let $0<\de<\de_n = (T^{\ast} - t_n)^{2\alpha}
\leq (T^*-t_n)$ (recall $2\alpha \geq 1$). Raising to the power
${\frac{(d+2)}{d}}$, dividing (\ref{limsup1}) by $l(J_n)^2$ and
integrating both sides with respect to $t$ on $(t_n, t_n + \de)$, we
obtain{\footnote{To be strictly correct, we should raise both sides
in (\ref{allR2Sup:1})-(\ref{allR2Sup}) to the power $\frac{d}{d+2}$
to obtain norms so that we can apply the triangle inequality.}}
\begin{equation}
 \label{allR2Sup:1}
\left(\frac{\ep}{2}\right)^{\frac{d+2}{d}}
\frac{\de}{(T^{\ast}-t_n)^{2\alpha}} \leq \int\limits_{(t_n, t_n +
\de)\times J_n} |P_{L(t_n)} u(x,t_n)|^{\frac{2(d+2)}{d}} \, dx \, dt
\end{equation}
\begin{equation}
 \label{allR2Sup:2}
 \leq \int\limits_{(t_n, t_n + \de)\times J_n}
|P_{L(t_n)}u(x,t_n) - P_{L(t)} u(x,t)|^{\frac{2(d+2)}{d}} \, dx \,
dt \,
\end{equation}
\begin{equation}
\label{allR2Sup} + \int\limits_{(t_n, t_n + \de)\times J_n}
|P_{L(t)} u(x,t)|^{\frac{2(d+2)}{d}}\, dx \, dt = I + II.
\end{equation}
Using the same estimate as in (\ref{E:difference}), we get
$$
|P_{L(t_n)}u(x,t_n) - P_{L(t)} u(x,t)| \leq 2 \, \Vert u_0
\Vert_{L^2(\cR^d)} \, [2 L(t_n+\de)]^{d/2}
$$
by H\"older's inequality and the conservation of mass.  Then the bound on term
$I$ in (\ref{allR2Sup}) is obtained by using the definition of
$L(t)$ and the bound on $l(J_n)$
$$
I \leq \int\limits_{(t_n, t_n + \de)\times J_n} \left(2 \,\Vert u_0
\Vert_{L^2(\cR^d)}\right)^{\frac{2(d+2)}{d}} \, [2 L(t_n +
\de)]^{d+2} \, dx \, dt
$$
$$
\leq (2 \,\Vert u_0 \Vert_{L^2(\cR^d)} )^{\frac{2(d+2)}{d}}
\frac{\kappa(\ep)^{d+2}}{(T^*-(t_n+\de))^{\alpha(d+2)}} \, |J_n| \,
\de
$$
$$
\leq \frac12 \, \left(\frac{\ep}{2}\right)^{\frac{d+2}{d}} \frac{\de
\,(T^{\ast}-t_n)^{d\alpha}}{(T^{\ast}-(t_n+\de))^{\alpha(d+2)}}.
$$
The second term in (\ref{allR2Sup}) is estimated by the space-time
norm on all $\cR^d$ ($J_n \subset \cR^d$)
$$
II \leq \int\limits_{(t_n, t_n + \de)\times \cR^d}
|u(t,x)|^{\frac{2(d+2)}{d}} \, dx \, dt = \Vert u
\Vert^{\frac{2(d+2)}{d}}_{L^{\frac{2(d+2)}{d}}{((t_n, t_n +
\de)\times \cR^d})}.
$$
Substituting all above estimates into (\ref{allR2Sup}), we obtain
\begin{align*}
\left(\frac{\ep}{2}\right)^{\frac{d+2}{d}}
\frac{\de}{(T^{\ast}-t_n)^{2\alpha}} &\leq \frac12 \,
\left(\frac{\ep}{2}\right)^{\frac{d+2}{d}} \frac{\de
\,(T^{\ast}-t_n)^{d\alpha}}{(T^{\ast}-(t_n+\de))^{(d+2)\alpha}}\\
& + \,\left[F(t_{n}+ \delta) - F(t_{n}) \right],
\end{align*}
which implies
\begin{align*}
\left(\frac{\ep}{2}\right)^{\frac{d+2}{d}} &\leq \frac12 \,
\left(\frac{\ep}{2}\right)^{\frac{d+2}{d}}
\frac{(T^{\ast}-t_n)^{(d+2)\alpha}}{(T^{\ast}-(t_n+\de))^{(d+2)\alpha}}\\
& + \, (T^{\ast}-t_n)^{2\alpha} \, \left( \frac{F(t_{n} + \delta) -
F(t_n)}{\delta}\right).
\end{align*}
Taking $\de \searrow 0$ ($n$ is fixed) and recalling that $F'$
exists at $t_n$, then absorbing the first term on the right into the
left-hand side,
we obtain
\begin{equation}
 \label{DerivativeGrowth}
\frac1{2} \left(\frac{\ep}2\right)^{\frac{d+2}{d}}
\frac{1}{(T^{\ast}-t_n)^{2\alpha}} \leq  F'(t_n),
\end{equation}
which gives (\ref{L4boundSup}).
\end{proof}

The preceding result shows that the time derivative of the Strichartz
norm is lower bounded along the sequence of times
where we have tight mass concentration. If we assume that the tight mass
concentration persists as in the conclusion of Corollary
\ref{CorThickTime}, we can integrate to obtain lower bounds on the
Strichartz norm itself.

\begin{lemma}\label{L:ExtLocal} 
Suppose that instead of the concentration (\ref{Mass1Sup}) with
$\limsup$, we have the concentration (\ref{refinedMass2proj}) with
the thickened time intervals $\{I_n\} \subset [t_n, t_{n+1})$  as in Corollary
\ref{CorThickTime}, i.e., for some $0 < \sigma < 1$ there exist
$\tilde\sigma > 0$ and cubes $\{J_n\} \subset \cR^d$ with $l(J_n) <
(T^* -t)^\alpha$ and
$|I_n| = \sigma (t_{n+1} - t_n)$ such that
\begin{equation}
 \label{refinedConc}
\lim_{n \to \infty} \inf_{t \in I_n} \int\limits_{J_n} |P_{L(t)} u(x, t)|^2 \, dx
\geq (1-\tilde\sigma)\, {\ep}.
\end{equation}
Then
\begin{equation}
 \label{ThickDerivative}
F'(t) \geq \frac{c(\sigma,\ep)}{(T^*-t)^{2\alpha}} \quad \text{for
a.e.} ~ t \in I_n.
\end{equation}
Furthermore, for all $\ds t \in \cup_{n} \,I_n$ we have
\begin{equation}
\label{LowerBound1} F(t) \gtrsim \left\{
\begin{array}{ll}
\ds \frac{1}{(T^*-t)^{2\alpha -1}} + \text{const}, &~\alpha > 1/2,\\
|\ln (T^*-t)| + \text{const}, &~ \alpha=1/2.
\end{array}
\right.
\end{equation}
\end{lemma}

\begin{proof}
Denote
$(t_n,t_n+\delta_n) = \text{int\,} I_n$, the interior of $I_n$. Take
$t \in (t_n, t_n+\delta_n)$ and repeat previous proof for this $t$
to obtain (\ref{ThickDerivative}) (note that the first step of
shifting $t_n$ to $t_n^*$ in order to have differentiability of $F$
available is not needed here, we may initially consider $t \in I_n$
such that $F'(t)$ exists).

For the second statement fix $t \in I_n$ and observe that since $F$
is increasing, we have $F(t) - F(t_n) \geq \int_{t_n}^{t} F'$.
Integrating the expression from (\ref{ThickDerivative}) (for $\alpha
> \frac12$) we obtain
$$
F(t) \geq \frac{c}{(T^* - t)^{2\alpha-1}} - \frac{c}{(T^* -
t_n)^{2\alpha-1}} + F(t_n),
$$
where $c=c(\alpha, \epsilon, \sigma)$. Next observe that $F(t_n)
\geq F(t_{n-1}+ \delta_{n-1}) \geq
\int_{t_{n-1}}^{t_{n-1}+\delta_{n-1}} F' + F(t_{n-1})$. Iterating
this process, say, till $t_{n-k} = T^*-1$, and using the property
that $\delta_{i} = \sigma(t_{i+1} - t_{i})$, we obtain that
$$
F(t) \geq \frac{c}{(T^* - t)^{2\alpha-1}} + \text{const}.
$$
Making appropriate changes in the integration for $\alpha =
\frac12$, we obtain the second part of (\ref{LowerBound1}).
\end{proof}
Observe that for $t \in (t_n, t_{n+1})\setminus I_n$, we have the
following estimate on $F$:
$$
F(t) \gtrsim \left\{
\begin{array}{ll}
\ds \frac{1}{(T^*-(t_n+\de_n))^{2\alpha -1}} + \text{const}, &~\alpha > 1/2,\\
|\ln (T^*-(t_n+\de_n))| + \text{const}, &~ \alpha=1/2,
\end{array}
\right.
$$
by using that $F$ is increasing on the compliment of $I_n$ in $(t_n,
t_{n+1})$.

\begin{corollary}
In analogy with Corollary \ref{GeneralCase} the statement of
Proposition \ref{propSup} can be generalized to include not only the
polynomial powers in the window of concentration in (\ref{Mass1Sup})
but a more general dependence on $(T^* - t)$. Suppose that both the
concentration window in (\ref{Mass1Sup}) is $l(J) < g(T^*-t)$ and
$\ds L(t) = \frac12 \frac{\kappa(\ep)}{g(T^*-t)}$, where the
function $g$ can be written as $g(T^*-t) = [-(\partial_t
G)(T^*-t)]^{-1/2}$ for some $C^1$-function $G$ with the properties
that as $t \to T^*$ both $G(T^*-t) \to \infty$ and $[-(\partial_t
G)(T^*-t)] \to \infty$.

Then the conclusion in (\ref{L4boundSup}) modifies to
\begin{equation}
 \label{GeneralBound}
F(t_n) \gtrsim G(T^*-t_n).
\end{equation}
\end{corollary}

For example, $g(T^* - t) = (T^* - t)^{\al} |\ln(T^*-t)|^{-\gamma}$
with $\alpha > 1/2$ and $\gamma \in \cR$, or $\alpha = 1/2$ and
$\gamma \geq 0$ \footnote{in fact, $\gamma > -1/2$ satisfies the
condition on $G$, however, for $-1/2 < \gamma < 0$ the window of
concentration is wider than parabolic} would satisfy the above
conditions. The last case produces the logarithmic divergence
$|\ln(T^*-t)|^{2\gamma +1}$ of the Strichartz norm $\ds \Vert u
\Vert_{L^{\frac{2(d+2)}{d}}{([0,t])\times \cR^d})}$.

\begin{proof}
To prove this general statement we repeat the proof of the above
Proposition with appropriate modifications and instead of
(\ref{DerivativeGrowth}) we arrive to
$$
F'(t_n) \gtrsim \frac1{[g(T^*-t_n)]^2}.
$$
Proceeding as in Lemmas \ref{L:ExtLocal}, and using the definition
of $g$, we obtain (\ref{GeneralBound}).

\end{proof}

\begin{corollary}
 \label{C:log}
If the blow up time $T^* < \infty$ for $NLS_{p}^{\pm}(\cR^d)$ with
$p=\frac4{d}+1$, then the diagonal Strichartz norm $\ds \Vert u
\Vert_{L^{\frac{2(d+2)}{d}}{([0,t])\times \cR^d})}$ explodes at
least as fast as $|\ln(T^* - t)|$.
\end{corollary}

The proof follows from \cite{B98}, \cite{BV05}, where it is shown
that finite time blowup solutions parabolically concentrate in
$L^2$, and the previous corollary with $g(T^*-t) = (T^*-t)^{1/2}$.

\end{document}